\documentclass[11pt,a4paper]{article}
\usepackage{amsmath,amssymb}
\usepackage{amsmath}

\usepackage{color}
\usepackage[colorlinks,linkcolor=blue,citecolor=red]{hyperref}

\DeclareMathSymbol{\bbbr}{\mathalpha}{AMSb}{"52}
\DeclareMathSymbol{\bbbc}{\mathalpha}{AMSb}{"52}

\newtheorem{theorem}{Theorem}

\newtheorem{corollary}[theorem]{Corollary}

\newtheorem{definition}{Definition}

\newtheorem{lemma}[theorem]{Lemma}

\newtheorem{proposition}[theorem]{Proposition}

\begin{document}

\title{Hexagonal geodesic 3-webs}

\author{{\Large Sergey I. Agafonov}\\
\\
Department of Mathematics,\\
S\~ao Paulo State University-UNESP,\\ S\~ao Jos\'e do Rio Preto, Brazil\\
e-mail: {\tt sergey.agafonov@gmail.com} }
\date{}
\maketitle
\unitlength=1mm

\vspace{1cm}

\begin{abstract}

\bigskip We prove that a surface carries a hexagonal 3-web of geodesics if and only if the geodesic flow on the surface admits a cubic first integral and show that the system of partial differential equations, governing  metrics on such surfaces, is integrable by generalized hodograph transform method. We present some new local examples of such metrics, discuss known ones, and establish an analogue of the celebrated Graf and Sauer Theorem for Darboux superintegrable metrics.     

\medskip


\medskip

\end{abstract}

%\vspace{-7mm}
%\newpage

%\tableofcontents

\section{Introduction}
Suppose that you are constructing a building of an original design: its roof is a smooth non-planar surface. 
You need a supporting framework, whose basic elements are bent timbers. For rigidity, you arrange them in a triangle grid. Economizing the material or/and  maximizing the construction strength towards the transverse loading, you conclude that the timbers have to follow geodesics of the roof surface. 
Which surface forms suit your roof? A bit surprising is that you are not completely free to choose! 

Let us phrase the question in the language of differential geometry: the framework becomes a triangle net on the surface, the framework elements form three discrete families of geodesics, intersecting in the net vertices. Letting the net mesh tend to zero, we obtain three geodesic foliations. 
A configuration of $d$ foliations is called {\it d-web}. Thus we have a 3-web of geodesics on our surface. The triangular combinatorics of the net survives in the disguise of {\it hexagonal} web. Topologically, this means that for each point $p$ on the surface any sufficiently small (but finite!) triangle, formed by geodesics of the foliations and  having the point $p$ as one of its vertices, can be completed to the curvilinear hexagon, whose sides are geodesics of the foliations and whose "large" diagonals are the tree foliation geodesics meeting at $p$. This non-trivial incidence relation amounts to vanishing of the Blaschke curvature of the 3-web. 

Thus, for the roof, we need a surface carrying a hexagonal geodesic 3-web. The history of studying such webs stretches back to the end of 19th century: Finsterwalder \cite{Fi-99} observed that any surface of revolution admits a one-parametric family of hexagonal geodesic 3-webs, and asked for more examples. His webs are constructed as follows: mark a point $p$ on the surface of revolution, draw the meridian $m$ and two geodesics, forming equal angles with $m$ at $p$, and rotate the surface around its axis to obtain 3 families of geodesics. The obtained 3-web is symmetric with respect to reflexion in any plane containing the axis. This immediately implies the hexagonality. 

The first example of geodesic hexagonal 3-web on a   surface without an infinitesiml isometry was published by St\"ackel \cite{St-03}.\footnote{St\"ackel communicated that the example was found by his student Ahl.} The surface turned out to be  a Lie spiral surface. 
Lie spiral surfaces are straightforward  generalizations of surfaces of revolution: instead of pure rotation of generating curves, one applies rotation combined with a homothety, whose center lies on the symmetry axis, at the rate $\exp(\alpha \theta)$, where $\alpha $ is constant and $\theta$ is the rotation angle. 

Further examples were provided by Sauer \cite{Sa-26}: he described surfaces of revolution, carrying hexagonal geodesic 3-webs of not Finsterwalder type and a class of hexagonal geodesic 3-webs on Lie spiral surfaces. Unfortunately, the treatment was incomplete and the author missed most of the webs on his surfaces of revolution and a class of Lie spiral surfaces. 
By clever ad hoc substitutions, Volk \cite{Vo-29} managed to reduce partial differential equations (PDEs), governing the metrics of searched for surfaces, to integrable cases of ordinary differential equations (ODEs). Most of the obtained examples were Liouville surfaces. 

Mayrhofer \cite{Ma-31} raised the issue of how many hexagonal geodesic 3-webs a surface can carry and obtained the sharp upper bound:  the family of such webs can depend on at most 9 parameters. The bound is realized only on the surfaces of constant Gaussian curvature.
 
The recent interest to hexagonal geodesic 3-webs is motivated by freeform architecture \cite{DPW-11,PHD-10}, as roughly outlined  in the first paragraph of this Introduction.   

It is remarkable, that the old problem of hexagonal geodesic 3-webs is closely related to another classical problem, namely that of integrability of the geodesic flow by first integrals, polynomial in momenta (or, equivalently, in velocity). Namely, in this paper we prove that a surface carries a hexagonal 3-web of geodesics if and only if the geodesic flow on the surface admits a cubic first integral. The relation of the web and the integral is quite natural: the web directions coincide with the zero directions of the integral. 

The problem of geodesic flow integrability by polynomial integrals dates from the 19th century classical works  of Dini and Darboux \cite{Di-69,Da-91}, where the local theory for the case of linear and quadratic integrals were developed (see also \cite{Ko-82} for a modern treatment and the review \cite{BMF-98} for more references). By the Noether Theorem, the linear integrals are in one-to-one correspondence with infinitesimal isometries, called also Killing vector fields of the surface. Therefore it is immediate that the dimension of the space of linear integrals is either 0 (generic case), or 1 (for the surfaces locally isometric to surfaces of revolution),   or 3 (for surfaces of constant Gaussian curvature).   The problem of dimension of the space of quadratic integrals was settled by Koenigs \cite{Ko-96}: the possible dimension is 1, 2, 3, 4, or 6.  

The case of cubic integrals is much more involved computationally. The space of cubic integrals is at most 10-dimensional, the bound realized on surfaces of constant Gaussian curvature. Moreover, it was shown in \cite{Kr-08} that for the case of non constant Gaussian curvature, the dimension is at most 7, and conjectured that the sharp bound is actually 4. The conjecture may be forced only by the lack of examples: all known examples support the claim. A classification of metrics admitting a Killing vector fields and at least one non-trivial cubic integral (i.e. non divisible by a linear one) was obtained in  \cite{MS-11} (see also \cite{VDS-15}). Under this hypothesis, the space of cubic integrals is indeed 4-dimensional.

There is another integrable system behind the problem of hexagonal geodesic 3-web, this time an infinite-dimensional one. If a surface carries a hexagonal geodesic 3-webs, the surface metric is subject to some system of partial differential equations of hydrodynamic type.  This system appeared as early as 1903 in a short note \cite{St-03} by St\"ackel, but its properties were not studied in detail neither by St\"ackel himself nor by other authors.       
We show that this system is diagonalizable and semi-Hamiltonian. Therefore it is integrable by generalized hodograph transform, discovered by Tsarev 
\cite{Ts-85}. 

If a given metric admits a polynomial integral of geodesic flow, the coefficients of the polynomial satisfy a linear overdetermined system of PDEs. If one considers the metric as unknown, then, for the coefficients and for the metric, one gets also a system of hydrodynamic type, which can be written in Cauchy form. In particular, this ensures a local existence of metrics, whose geodesic flow has a polynomial integral of any prescribed degree $n$ (see \cite{Te-97}). It is interesting that this system for metric and coefficients is also diagonalizable and semi-Hamiltonian (see \cite{BM-11,MP-17}). 

In this paper, we also study symmetry properties of the system for metric. Apart from the discrete permutation symmetry $S_3$, caused by ambiguity in the choice of local coordinates, adjusted to the web, the system admits a 4-dimensional algebra of infinitesimal symmetries. Acting on solutions, each symmetry  induces either isometry or homothety on the corresponding surface.   We provide a classification of the one-dimensional subalgebras and clarify the geometrical meaning of solutions, invariant with respect to the operators, generating homotheties: these solutions define metrics on the surfaces, locally isometric to the Lie spiral surfaces.

Infinitesimal projective transformations, or projective vector fields, are such vector fields on a surface, whose flow moves around the geodesics as  unparametrized curves. Lie showed \cite{Li-82} that the dimension of the projective algebra $ \mathfrak{p}(g)$ is 0, 1, 2, 3 or 8, where only the surfaces of constant Gaussian curvature have  $\dim \mathfrak{p}(g)=8$. If $\dim \mathfrak{p}(g)\ge 2$ then the metric necessarily admits a Killing vector field (see \cite{BMM-08}) and, consequently, at least a 1-parameter family of hexagonal geodesic 3-webs, whereas for  $\dim \mathfrak{p}(g) =3$ the metric is {\it Darboux superintegrable} (i.e. the space of quadratic first integrals is 4-dimensional) and the family of hexagonal geodesic 3-webs is 3-parametric. We describe such webs for $\dim \mathfrak{p}(g) = 2,3$ and, for Darboux superintegrable metrics, establish an exact analogue of the celebrated Theorem of Graf and Sauer (see \cite{GS-24}). Recall that the Graf and Sauer Theorem claims that any 3 concurrent lines $l_1,l_2,l_3$, of a linear hexagonal 3-web on the projective plane $\mathbb{P}^2$ are dual to the points, where the line, dual to the point $p=l_1\cap l_2 \cap l_3$, meets some fixed (possibly singular) cubic in the dual plane $(\mathbb{P}^2)^*$. For surfaces with $\dim \mathfrak{p}(g) =3$, we represent  the dual space (i.e. the space of geodesics) by some quadric $Q$ in $\mathbb{P}^3$. Then the geodesic foliations correspond to curves on $Q$, the pencil of concurrent geodesics with the vertex $p$ being represented by a plane conic $C_p$ on $Q$.   The description of hexagonal geodesic 3-webs in dual terms is as follows: fix two planes $P_0$ (this one is the same for all webs) and $P_1$ (which varies with the web) in $\mathbb{P}^3$, cutting the quadric $Q$ in two conics $C_0$ and $C_1$, then the 2 duals of 3 concurrent geodesics of the web, passing through a point $p$, are the points, where the plane of $C_p$ meets $C_1$, and the third dual is the point, where the plane of $C_p$ touches $C_0$.

Finally, we present some new local examples of hexagonal geodesic 3-webs, obtained via group-invariant solutions of the PDE system for metric, or via simple wave solutions. In particular, we rectify some claims made in \cite{Sa-26}, present a class of Lie spiral surfaces, missed in \cite{Sa-26}, and give new forms for some known examples.

\section{Integrable quasilinear system behind hexagonal geodesic 3-webs}
Suppose that a surface $(M^2,g)$ carries a hexagonal 3-web ${\mathcal G}_3$ formed by geodesic foliations ${\mathcal F}_i$, $i=1,2,3$.  In what follows, we  suppose ${\mathcal G}_3$ to be non-singular at a generic point $p_0\in M^2$, which means that foliations are pairwise transverse. Then the web hexagonality implies the existence of locally defined non-constant smooth functions $u,v,w$ such that:\\
1) $u,v,w$ are constant along the leaves of ${\cal F}_1, {\cal F}_2$ and ${\cal F}_3$ respectively,\\
2) $u+v+w\equiv 0$ (see \cite{BB-38}).

 Since the web is non-singular at $p_0$, any two of these functions serve as local coordinates on $M^2$. Let us write the metric $g$ in the coordinates $u,v$: 
$$
ds^2=E(u,v)du^2+2F(u,v)dudv+G(u,v)dv^2.$$
The web leaves are solutions of ordinary differential equations 
$$
du=0, \ \ \ dv=0,\ \ \  dv=-du.
$$ Taking into account the equation for unparametrized geodesics
\begin{equation}\label{geoEQ}
\frac{d^2v}{du^2}=-\Gamma ^2_{11}+(\Gamma ^1_{11}-2\Gamma ^2_{12})\frac{dv}{du}-(\Gamma ^2_{22}-2\Gamma ^1_{12})\left( \frac{dv}{du}\right)^2+\Gamma ^1_{22}\left( \frac{dv}{du}\right)^3,
\end{equation}
where $\Gamma ^i_{jk}$ are Christoffel's symbols of the Levi-Civita connection,
one gets
$$
\Gamma ^2_{11}=\Gamma ^1_{22}=\Gamma ^1_{11}-2\Gamma ^2_{12}+\Gamma ^2_{22}-2\Gamma ^1_{12}=0.
$$
Invoking the formulas for the Christoffel symbols, we arrive at the following system of PDEs:
\begin{equation}\label{hydro}
\begin{array}{l}
2EF_u-FE_u-EE_v=0\\
\\
2GF_v-FG_v-GG_u=0\\
\\
GE_u+EG_v-2F(F_u+F_v)+(3F-2G)E_v+(3F-2E)G_u=0.
\end{array}
\end{equation}
Observe that system (\ref{hydro}) is quasilinear and belongs to the so-called hydrodynamic type, i.e. for ${\bf Z}:=(E,F,G)^T$, it has the form
\begin{equation}\label{hydrotype}
A({\bf Z}){\bf Z}_u+B({\bf Z}){\bf Z}_v=0,
\end{equation}
where $A,B$ are $3\times 3$ matrices depending only on ${\bf Z}$.
Let as recall the basic definitions and concepts of the theory of such systems.
\begin{definition}
System (\ref{hydrotype}) is said to be  {\it diagonalizable} if one can choose dependent variables $R^i$, called {\it Riemann invariants}, so that the system assumes the form
\begin{equation}\label{n}
R^i_u=\lambda^i(R) R^i_v,
\end{equation}
$i=1, 2, 3$ (no summation over $i$).
\end{definition}
This type systems govern a wide range of problems in pure and applied mathematics.

\begin{definition}
System (\ref{n}) is called {\it semi-Hamiltonian} if its characteristic speeds $\lambda^i$ satisfy the   constraints
\begin{equation}\label{semiH}
\partial_k \left(\frac{\partial_j\lambda
^i}{\lambda^j-\lambda^i}\right)=\partial_j \left(\frac{\partial_k\lambda
^i}{\lambda^k-\lambda^i}\right),
\end{equation}
here $\partial_i=\partial/\partial {R^i}$.
\end{definition}
Tsarev \cite{Ts-85} showed  that this property implies integrability:  semi-Hamiltonian systems (\ref{n}) possess infinitely many conservation laws and commuting flows, and can be solved by the generalised hodograph method.

\begin{theorem}\label{Tsarevintegrability}
System (\ref{hydro}) is diagonalizable and semi-Hamiltonian.
\end{theorem}
{\it Proof:} The characteristic speeds $\lambda^i$ satisfy the characteristic equation $\det[B(z)+\lambda A(z)]=0$, i.e.
\begin{equation}\label{char}
G\lambda ^3+(F-2G)\lambda^2+(F-2E)\lambda+E=0.
\end{equation}
Riemann invariants $R^i$ are first integrals of the 2-dimensional distributions generated by two vector fields
$\xi_j,\xi_k$, where $j,k$ are distinct and $\xi_l$ verifies $[B(z)+\lambda^l A(z)]\xi_l=0$ (no summation over $l$).  We choose the following rational function for $R^3$:
\begin{equation}\label{R3}
\begin{array}{l}

R^3= {\scriptstyle
{ [2  \lambda^1  \lambda^2+ {(\lambda^2)}^2- \lambda^1-2  \lambda^2]^2 [2  \lambda^1  \lambda^2+ {(\lambda^1)}^2-2  \lambda^1- \lambda^2]^2 F}} /
 \\
 \ \ \ \ \left( \right. {\scriptstyle
 [2  \lambda^1  \lambda^2- \lambda^1- \lambda^2-1] [ \lambda^1  \lambda^2+ \lambda^1+ \lambda^2-2] [ \lambda^1  \lambda^2-2  \lambda^1-2  \lambda^2+1] [2 \lambda^1 \lambda^2-\lambda^1-\lambda^2+2]}\\
\left. {\scriptstyle
 \ \ \ \  [2 {(\lambda^1)}^2 \lambda^2+2 \lambda^1 {(\lambda^2)}^2+2 (\lambda^1+\lambda^2)
 -5 \lambda^1 \lambda^2-{(\lambda^1)}^2-{(\lambda^2)}^2]
 }\right),
 \end{array}
\end{equation}
and  obtain $R^1,R^2$ from the formula for $R^3$ by cyclic permutations of $\lambda$-s. Now the semi-
Hamiltonian property (\ref{semiH}) is checked by direct computation.
\hfill $\Box$\\

\smallskip
\noindent {\bf Remark 1}. The characteristic speeds $\lambda^1,\lambda^2,\lambda^3$ are not independent functions on the hodograph space, they are subject to the following identity
\begin{equation}\label{lambdas}
\lambda^1\lambda^2+\lambda^2\lambda^3+\lambda^3\lambda^1+
\lambda^1+\lambda^2+\lambda^3=2+2\lambda^1\lambda^2\lambda^3.
\end{equation}
Condition (\ref{semiH}) was checked with the help of symbolic computation software, namely Maple.

\medskip
\noindent {\bf Remark 2}.
System (\ref{n}) is {\it linearly degenerate} if its characteristic speeds verify the conditions
$$
\partial_i\lambda^i=0,
$$
no summation,  $i=1, 2, 3$. Linear degeneracy prevents the breakdown of smooth initial data, which is typical for genuinely non-linear systems of type (\ref{n}) (see e.g. \cite{RS-67}). System (\ref{hydro}) is not linearly degenerate, therefore one expects that its generic solution defines a smooth surface, carrying hexagonal geodesic 3-web, only locally.

\section{Symmetry properties}\label{symsec}
System (\ref{hydro}) has obvious discrete symmetries related to the ambiguity in the choice of the web first integrals $u,v,w$ as local coordinates. The following substitutions generate the group of discrete symmetries of (\ref{hydro}): \\

transposition   $T_{uv}:\ \bar{u}=v,\ \bar{v}=u, \ \bar{E}=G, \ \bar{G}=E, \ \bar{F}=F$,\\

substitution    $T_{uw}:\ \bar{u}=w=-u-v,\ \bar{v}=u, \ \bar{E}=G, \ \bar{G}=E+G-2F, \ \bar{F}=G-F$,\\ 

changing the sign $\ \bar{u}=-u,\ \bar{v}=-v, \ \bar{E}=E, \ \bar{G}=G, \ \bar{F}=F$. \\

The system also has a non-trivial Lie point symmetry algebra.  Let us recall the basic concepts of the Lie theory and describe this algebra.

Any solution to (\ref{hydro}) defines, at  least locally, some surface $\Sigma$ parametrized by $u,v$:
\begin{equation}\label{solsurf}
\sigma : V\subset \mathbb R^2 \to \Sigma \subset \mathbb R^5, \  \ \ \ (u,v) \mapsto (u,v,E(u,v),F(u,v),G(u,v)).
\end{equation}
Consider a local flow of some vector field $X$ on $\mathbb R^5$
\begin{equation}\label{infinitesimal}
X=\alpha\partial_u+\beta\partial_v+\varepsilon \partial_E+\varphi \partial_F+\gamma\partial_G.
\end{equation}
This flow deforms a surface $\Sigma $, defined  via (\ref{solsurf}) by a solution to (\ref{hydro}), to some  surface $\tilde{\Sigma}$. For the flow transforms, close to the identity, the surface $\tilde{\Sigma}$ remains a graph of some map 
$$
\tilde{\bf Z} : V\subset \mathbb R^2 \to \mathbb R^3, \ \ \ (u,v) \mapsto (\tilde{E}(u,v), \tilde{F} (u,v),\tilde{G} (u,v))
$$ 
and therefore can be parametrized by a suitable $\tilde{\sigma}$ of the form (\ref{solsurf}).

\begin{definition}
The vector field (\ref{infinitesimal}) is an infinitesimal symmetry of system (\ref{hydro}), if its local flow moves a surface $\Sigma $, defined by any solution to (\ref{hydro}), to a surface $\tilde{\Sigma}$ that again is defined by some solution to (\ref{hydro}).
\end{definition}
One easily checks that the following four vector fields are infinitesimal symmetries of  (\ref{hydro}):
\begin{equation}\label{sym}
\{T_1=\partial_u,\ T_2=\partial_v,\ D_1=u\partial_u+v\partial_v,\ D_2=E\partial_E+F\partial_F+G\partial_G \}
\end{equation}
It turns out that this list exhausts all the point symmetries of system (\ref{hydro}).

\begin{theorem}
 The symmetry algebra of PDE system (\ref{hydro}) is generated by the vector fields (\ref{sym}).
\end{theorem}
{\it Proof:}  The vector field (\ref{infinitesimal}) is an infinitesimal symmetry of (\ref{hydro}) if and only if it satisfies a system of so called {\it determining equations}. This PDE system is linear and strongly overdetermined, a straightforward analysis of its compatibility conditions yields its general solution as a linear combination of the fields (\ref{sym}). This analysis is considerably simplified by the following observation: the flow leaves  invariant our hexagonal geodesic 3-web formed by integral curves of ODEs $du=dv=du+dv=0$. Therefore the truncated operator $\alpha\partial_u+\beta\partial_v$ is a linear combination of $T_1,$ $T_2,$ and $D_1$.  (See \cite{Ol-93} for modern exposition of the Lie theory and for examples of solving determining equations). 
\hfill $\Box$
\begin{definition} A map
$$
\varphi: (M^2,g)\to (\tilde{M}^2,\tilde{g}),
$$ 
is called a projective map of the surface $(M^2,g)$ to the surface $(\tilde{M}^2,\tilde{g})$ if 
it sends unparameterized geodesics of $(M^2,g)$ to unparameterized geodesics of $(\tilde{M}^2,\tilde{g})$.\\
A projective map is called a homothety if $\varphi ^*\tilde{g}=const\cdot g$ 
\end{definition}
The symmetry algebra (\ref{sym}) generates a pseudogroup of diffeomorphisms, each element of this pseudogroup moving around solutions of (\ref{hydro}). If $(M^2,g)$ is a surface defined by a solution of (\ref{hydro}) and $X$ an operator of the symmetry algebra, then for sufficiently small $t\in \mathbb R$, the  exponential $\varphi_t=\exp (tX)$ induces a  map
\begin{equation}\label{syminduced}
\tilde{\varphi}_t: (M^2,g)\to (M^2_t,g_t)
\end{equation}  
 to some surface $(M^2_t,g_t)$,
 also carrying a hexagonal geodesic 3-web. 
\begin{theorem}\label{projsym}
Let $(M^2,g)$ be a surface with a metric $g$ defined by some solution of (\ref{hydro}). Then
for any vector field  $X$ of the symmetry algebra (\ref{sym}), the map (\ref{syminduced}) induced by the exponential $\varphi_t=\exp (tX)$ is projective. Moreover, if this map is not an isometry then it is a homothety. 
\end{theorem}
{\it Proof:} Consider a geodesic on $(M^2,g)$. This geodesic, as a curve in the $u,v$-coordinate chart, is a solution to the second order ODE (\ref{geoEQ}), which takes the form 
\begin{equation}\label{geoEQ2}
\frac{d^2v}{du^2}=K(u,v)\frac{dv}{du}\left(1+\frac{dv}{du}\right),
\end{equation}
where 
$$
K(u,v)=\Gamma ^1_{11}-2\Gamma ^2_{12}=-\Gamma ^2_{22}+2\Gamma ^1_{12}=\frac{GE_u+3FE_v-2FF_u-2EG_u}{2(EG-F^2)}.
$$
Let us lift a point $(u,v)\in M^2$  to the parametrized surface (\ref{solsurf}), move this surface by $\varphi_t=\exp (tX)$ and then project the image back to the $uv$-plane. This defines the local map $\tilde{\varphi}_t$. Now it is straightforward that the map  $\tilde{\varphi}_t$, prolonged to the derivatives $\frac{dv}{du}$, $\frac{d^2v}{du^2}$, sends equation (\ref{geoEQ2}) to its counterpart for $(M^2_t,g_t)$. This proves the first claim of the Theorem. 

Now observe that $\tilde{\varphi}_t$, corresponding to the translations $T_1,T_2$ (see (\ref{sym})), acts trivially on the metric $g$, while the dilatations $D_1,D_2$ expand it by a constant factor. Therefore  $\tilde{\varphi}_t$ is either isometry or homothety.    
\hfill $\Box$

\smallskip

A vector field on the surface $(M^2,g)$, whose local flow maps unparameterized geodesics to unparameterized geodesics, is called {\it projective }. All projective vector fields on the surface $(M^2,g)$ form the projective symmetry algebra $\mathfrak{p}(g)$. 

It can happen that $\exp (tX)$ preserves the parametrized surface (\ref{solsurf}). Then the corresponding solution to (\ref{hydro}) is called {\it invariant} with respect to the operator $X$, and the diffeomorphism $\tilde{\varphi}_t$ is a projective transformation of $(M^2,g)$. Moreover, the vector field $X$ is tangent to the surface (\ref{solsurf}) and defines a {\it projective } vector field on $(M^2,g)$, which we denote by the same symbol $X$. 
       
It is natural to consider metrics, admitting hexagonal geodesic 3-webs, up to action of the pseudogroup generated by the symmetry algebra (\ref{sym}). In particular, we will consider invariant solutions up to this action. 
If two solutions to  (\ref{hydro}) are related by a transform from the above pseudogroup and one of the solutions is invariant with respect to an operator $X$, then the second solution is also invariant with respect to some $\tilde{X}$, where the one-dimensional algebras, generated by $X$ and $\tilde{X}$ are related by an inner automorphism of the symmetry algebra (see \cite{Ov-82}).  
\begin{proposition}\label{Liesysm}
If a solution to  (\ref{hydro}) is invariant with respect to a one-dimensional subalgebra of algebra (\ref{sym}), then, permuting $u$ and $v$ if necessary, one finds a constant $\kappa$ and an inner automorphism of symmetry algebra, sending this subalgebra to one, generated by an operator from the following list:
\begin{equation}\label{normalLiesysm}
\begin{array}{ll}
1) & \partial_u+\kappa \partial_v,\\
2) & \partial_u+\kappa \partial_v+E\partial_E+F\partial_F+G\partial_G,\\
3) & u\partial_u+v\partial_v+\kappa(E\partial_E+F\partial_F+G\partial_G).
\end{array}
\end{equation}
\end{proposition} 
{\it Proof:} With a small abuse of notation, we will denote by $\exp (tX)$ also the exponential of the inner automorphism generated by an element $X$.  Let $c_i$ be constants and $X=c_1T_1+c_2T_2+c_3D_1+c_4D_2$ be a symmetry operator. 

If $c_3=0$, then acting on $X$ by $\exp (tD_1)$ one can rescale $c_1$ and $c_2$ by the same factor. Note that if our solution to (\ref{hydro}) is invariant with respect to the operator $X$ then $X\wedge D_2\ne 0$ and at least one of the coefficients $c_1,c_2$ does not vanish. Suppose that $c_1\ne 0$. If $c_4=0$ then our operator can be rescaled to the form 1). If $c_4\ne 0$ then applying $\exp (tD_1)$ and multiplying $X$ by a  suitable constant factor, one arrives to the form 2). 

If $c_3\ne 0$ then applying $\exp (tT_1)$ and $\exp (tT_2)$ one kills coefficients $c_1$ and $c_2$ and gets the form 3) after a suitable rescaling. 
\hfill $\Box$\\ 
It is well known that a surface, whose metric admits a Killing vector field, can be locally immersed in Euclidean space $\mathbb R^3$ as a surface of revolution.  

\begin{lemma}\label{spiralLie}
Let $Y$ be a vector field on a surface such that the Lie derivative of its metric $g$ verifies the condition   $\mathcal{L}_Y(g)=g$. Then the surface can be locally immersed in Euclidean space $\mathbb R^3$ as a  Lie spiral surface. 
\end{lemma}
{\it Proof:} Choosing the axis of revolution to be the z-axis and taking the origin as the homothety center, one can parametrize a Lie spiral surface by $\theta, r$ as follows: 
$$
x=e^\theta r \cos (\alpha \theta),\ \ \ y=e^\theta r \sin (\alpha \theta), \ \ \ z=e^\theta W(r).
$$
The induced metric is 
$$
ds^2=e^{2\theta}\{[r^2(1+\alpha ^2)+W^2]d\theta^2+2[r+WW']d\theta dr+[1+(W')^2]dr^2\}. 
$$
On the given surface, let us choose local coordinates $u,v$ so that $Y=\partial_u.$
Then its metric has the form 
$$
d\bar{s}^2=e^u(h(v)du^2+2f(v)dudv+g(v)dv^2).
$$ 
There are (locally defined) functions $U(r),V(r),W(r)$ such that the coordinate transform 
$$
u=2\theta+U(r), \ \ \ v=V(r)
$$
brings the quadratic form $d\bar{s}^2$ to the quadratic form $ds^2$, 
if these function satisfy the following equations:
$$
\begin{array}{l}
4e^Uh(V)=(\alpha ^2+1)(W^2+r^2),\\
\\
2e^U[h(V)U'+f(V)V']=WW'+r,\\
\\
e^U[h(V)(U')^2+2f(V)U'V'+g(V)(V')^2]=(W')^2+1.
\end{array}
$$ 
This system is locally solvable for $U,V,W$. To check this, let us differentiate the first equation, solve it for $W'$, substitute thus found $W'$ into the second equation and solve it for $U'$. Now the third equation gives a quadratic equation for $V'$. Its discriminant $D$
$$
D=\frac{16h^2\alpha^6e^{2U}}{1+\alpha^2}([1+\alpha^2](f^2-gh)+gh-2fh'+(h')^2)([1+\alpha^2]r^2-4e^UhV)
$$ 
could be made non-negative by an appropriate choice of initial value for $U$. (In the above formula, the argument of $f,g,h,h'$ is $V$.)  
\hfill $\Box$
\begin{corollary}
If a solution to (\ref{hydro}) is invariant with respect to the one-dimensional subalgebra (\ref{normalLiesysm},2)) (with any $\kappa$) or (\ref{normalLiesysm},3)) with $\kappa\ne -2$, then the surface  with the corresponding metric can be realized locally as a  Lie spiral surface. 
\end{corollary}
In fact, for the case of operator (\ref{normalLiesysm},3)) holds true $\mathcal{L}_Y(g)=(\kappa+2)g$, where 
$Y=u\partial_u+v\partial_v$. 
 
\section{Integrability of geodesic flows via cubic first integrals}
Consider the geodesic flow on $(M^2,g)$ with the metric $g(du,dv)=Edu^2+2Fdudv+Gdv^2$. Let $(p,q)$ be the momentum so that $p$ is conjugate to $u$ and $q$ to $v$ respectively. Then the Hamiltonian reads as
$$
H =\frac{Gp^2-2Fpq+Eq^2}{2(EG-F^2)}.
$$
Vector fields $\xi\partial_u+\eta\partial _v$, tangent to the leaves of  a hexagonal geodesic 3-web on the surface  $(M^2,g)$, where $(E,F,G)$ verify (\ref{hydro}), are annihilators of the binary cubic form $du dv (du+dv)$:
\begin{equation}\label{webcubicform}
\xi\eta(\xi+\eta)=0.
\end{equation}
Under the canonical isomorphism $\psi : TM\to T^*M$, to the momentum $(p,q)$, there corresponds the vector $(\xi,\eta)=\frac{1}{EG-F^2}(Gp-Fq,Eq-Fp)$. Therefore equation (\ref{webcubicform}) can be rewritten in terms of $(p,q)$ as vanishing of the following polynomial, cubic in $(p,q)$:
\begin{equation}\label{webcubicint}
I=\mu \cdot(Gp-Fq)(Eq-Fp)[(G-F)p+(E-F)q],
\end{equation}
where $\mu$ is any non-vanishing function of $u,v$.
\begin{lemma}\label{LMcub}
Let $(E,F,G)$ be a solution to (\ref{hydro}). Then there is an "integrating" factor  $\mu=\frac{1}{EG-F^2}$ such that
the cubic polynomial (\ref{webcubicint}) is a first integral of the geodesic flow of the metric $g$. The factor  $\mu$ is defined up to multiplication by constant.\footnote{The author thanks the referee for indicating the explicit form of the integrating factor.}
\end{lemma}
{\it Proof:} The commuting condition $\{I,H\}=0$, being a homogeneous polynomial equation of fourth degree in $p,q$, splits to give 5 equations. Only two of them are independent. Resolving them for $\frac{\mu_u}{\mu},\frac{\mu_v}{\mu}$ and substituting for $E_u,F_u,G_u$ their expressions via $E_v,F_v,G_v$ obtained from (\ref{hydro}), one 
easily  checks that $d\ln (\mu)=-2d\ln(EG-F^2)$
on any solution to (\ref{hydro}).
\hfill $\Box$

\bigskip
\noindent It is remarkable that the above Lemma is invertible.
\begin{definition}
 We call a direction $[p:q]\in {\mathbb R\mathbb P^1}$ the real root of a homogeneous polynomial $I\in \mathbb R [p,q]$, if $I(p,q)=0$.
\end{definition}
Any cubic integral of geodesic flow is, by definition, homogeneous in $p,q$, the coefficients being smooth functions of a point $(u,v)$ on the surface.
\begin{theorem}\label{cubandflat}
A surface $(M^2,g)$ carries a hexagonal geodesic 3-web if and only if the geodesic flow on this surface admits a cubic first integral with 3 distinct real roots.
\end{theorem}
{\it Proof:} Due to Lemma \ref{LMcub}, it is enough to show that the existence of a cubic first integral of the geodesic flow implies the existence of a hexagonal geodesic 3-web.

A cubic integral determines, via the canonical isomorphism $\psi : TM\to T^*M$, a cubic binary form $\tau$ on $M$. Distinct real roots of the integral are mapped into 3 direction fields on $M$. The invariance of the integral along the geodesic flow implies that the 3-web, formed by integral curves of these 3 direction fields, is geodesic.

To calculate the Blaschke curvature of the web, let us choose local coordinates $(u,v)$ on $M$ so that $\tau$ is divisible by $du \cdot dv$. In these coordinates, the first integral takes the following form
$$
I=(Gp-Fq)(Eq-Fp)[L(Gp-Fq)+K(Eq-Fp)],
$$
where $L,K$ are functions of $u,v$,
and the web is formed by coordinate lines $u=const$, $v=const$ and by the integral curves of the ODE $Ldu+Kdv=0$. The Chern connection form  of the web is
$\Gamma =\frac{K_u}{K}du+\frac{L_v}{L}dv$.

 The commuting condition $\{I,H\}=0$ again splits into 5 equations. Their differential consequences imply vanishing of the Blaschke curvature form  $K_B=d(\Gamma)=0$. Thus the constructed geodesic 3-web is hexagonal. One possible way to check this implication is presented in Appendix.
\hfill $\Box$\begin{corollary}
If the space of cubic integrals of the geodesic flow on $(M^2,g)$ is $d$-dimensional then the family of hexagonal geodesic 3-webs on $M$  depends on $d-1$ essential parameters. In particular, if the set of hexagonal geodesic 3-webs on $M$ is finite then it contains only one element. 
\end{corollary}

\noindent {\bf Remark}. The existence of one-parameter family of hexagonal geodesic 3-webs on surfaces admitting a Killing vector field, indicated by Finsterwalder for surfaces of revolution, can be explained in a purely algebraic manner.  It is well known that the existence of Killing vector field implies the existence of linear integral (Noether Theorem), square of which gives a quadratic integral.  The metric $g$ provides another quadratic integral. Multiplying the quadratic integrals again by the linear one, we obtain a pencil of cubic integrals.

\section{Hexagonal geodesic 3-webs on surfaces with transitive projective symmetry algebra} 

Lie \cite{Li-82} classified projective algebras $\mathfrak{p}(g)$ on surfaces $(M^2,g)$. The possible dimensions of such non-trivial algebras are 1,2,3 or 8, the last case being that of constant Gaussian curvature. In this section we describe hexagonal geodesic 3-webs on the surfaces with $\dim \mathfrak{p}(g)=2,3$. We will use the slightly modified normal forms of projective connections obtained in \cite{BMM-08}.

\subsection{Hexagonal geodesic 3-webs on the surfaces with $\dim \mathfrak{p}(g)=3$}
For this dimension of projective symmetry algebra, one can choose local coordinates $x,y$ on the surface so that the equation of geodesics is
$$
\frac{d^2y}{dx^2}=\frac{1}{2}\frac{dy}{dx}+Je^{-2x}\left(\frac{dy}{dx}\right)^3,\ \ \ J={\rm const} \ne 0,
$$
with the following basis of the projective symmetry algebra
$$
\{\partial_y,\ \partial_x+y\partial_y,\ 2y\partial_x+y^2\partial_y \}
$$ (see \cite{BMM-08}).
Instead of the coordinate $x$ we prefer to use a new coordinate $z$, defined by the differential relation $\frac{dz}{z}=\frac{dx}{2}$. By adjusting the integration constant,
we obtain  the following equation of geodesics
\begin{equation}\label{eqgeo3}
z^3\frac{d^2y}{dz^2}=\epsilon\left(\frac{dy}{dz}\right)^3,\ \ \ \epsilon=\pm 1
\end{equation}
and the symmetry algebra generators
\begin{equation}\label{pg3}
\{\partial_y,\ z\partial_z+2y\partial_y,\ zy\partial_z+y^2\partial_y \}
\end{equation}
in the coordinates $z,y$.

Observe that the first two generators form a basis of a 2-dimensional subalgebra. According to the classical results of Lie, the second order equation (\ref{eqgeo3}) is integrable in quadratures (see \cite{Ol-93} for the modern exposition). In our case, the integration is straightforward: equation (\ref{eqgeo3}) is homogeneous ODE of the first order for $\frac{dy}{dz}$, which is integrable in closed form, namely one gets
$$
\left(\frac{dz}{dy}\right)^2= \frac{\epsilon}{z^2}+const.
$$
Integrating again one obtains the general solution in an implicit form
\begin{equation}\label{solgeo3}
k^2(y-l)^2-kz^2=\epsilon,
\end{equation}
where $k,l$ are integration constants. Note that special solutions $y=const$, which will play an important role in what follows, are obtained by limit $k\to \infty$.
Thus, the geodesics form a 2-parametric family of conics, the special solutions being considered as double lines $(y-l)^2=0$.
Let us consider the closure of this family, which we will call the {\it dual} space of $(M^2,g)$.
Any geodesic is a conic of the form
\begin{equation}\label{incidence}
Ay^2+2By+C+Dz^2=0.
\end{equation}
Taking into account the parametrization of the family by the parameters $k,l$ via (\ref{solgeo3}), one obtains immediately that the dual space is the quadric
\begin{equation}\label{quadric}
AC-B^2+\epsilon D^2=0.
\end{equation}
in $\mathbb R\mathbb P^3$ (a hyperboloid  of one sheet for $\epsilon=1$ and a hyperboloid  of two sheets for $\epsilon=-1$).
The special solutions correspond to the section of the above quadric by the plane $D=0$. This section is a smooth conic $c_0$.

To complete the duality picture, consider a point $(z,y)\in M^2$ on the surface and the pencil of geodesics centred at this point. To this pencil there corresponds the section of the quadric (\ref{quadric}) by the plane (\ref{incidence}). This plane intersects the plane $D=0$ along a line tangent to the conic $c_0$.

Now for any geodesic 3-web ${\cal{G}}_3$ on $M^2$ we have 3  arcs $\gamma_i, \ i=1,2,3,$ on the quadric  (\ref{quadric}). We call these arcs $\gamma_i$ the {\it dual focal curves} of the web ${\cal{G}}_3$. Any plane (\ref{incidence}) cuts these arcs in 3 points, corresponding to  the  3 web leaves  passing through $(z,y)$.
Hexagonal geodesic 3-webs enjoy a description similar to the one that Graf and Sauer provided for linear hexagonal 3-webs in the plane \cite{GS-24}.

\begin{theorem}\label{analogGS}
Suppose that $\dim \mathfrak{p}(g)=3$ and the coordinates $(z,y)$ on the surface $(M^2,g)$ are chosen as above. Then a geodesic 3-web ${\cal{G}}_3$ on $M^2$ is hexagonal if and only if one of its focal curves is an arc of the conic  $c_0$ and the other two are arcs of the section of the quadric (\ref{quadric}) by some fixed plane, different from the plane $D=0$.
\end{theorem}
{\it Proof:} The foliations of a geodesic 3-web on $M^2$ are  formed  by integral curves of  3   vector fields
\begin{equation}\label{webVfields}
\partial_z+P(z,y)\partial_y,\ \ \ \ \ \partial_z+Q(z,y)\partial_y,\ \ \ \ \ \partial_z+R(z,y)\partial_y,
\end{equation}
 satisfying the following system of decoupled equations:
\begin{equation}\label{3dimPQR}
z^3(P_z+PP_y)=\epsilon P^3,\ \ \ \ \ \ z^3(Q_z+QQ_y)=\epsilon Q^3,\ \ \ \ \ \ z^3(R_z+RR_y)=\epsilon R^3.
\end{equation}
The web is hexagonal if and only if its Blaschke curvature vanishes (see \cite{BB-38} for formulas). This constraint takes the form
\begin{equation}\label{cur03}
Pyy+Qyy+Ryy-\frac{3\epsilon}{z^3}(PPy+QQy+RRy)-\frac{3\epsilon}{z^4}(P+Q+R)+\frac{1}{z^6}(P^3+Q^3+R^3)=0.
\end{equation}
Analysis of compatibility of (\ref{cur03}) with system (\ref{3dimPQR}) implies that one of the slopes $P,Q,R$ vanishes identically. We may suppose that $R=0$, then we obtain a completely integrable Pfaff system for $P,Q,Py$:
\begin{equation}\label{3webeq3}
\begin{array}{l}
dP=\left(\frac{\epsilon P^3}{z^3}-PPy\right)dx+Pydy\\
\\
dQ=\left(\frac{Q^2Py}{P}-\frac{Q(\epsilon P^2Q+z^2P+z^2Q)}{z^3P} \right)dx +
\left( \frac{(P+Q)(\epsilon PQ+z^2)}{z^3P} -\frac{QPy}{P} \right)dy\\
\\
dPy=\left(\frac{Py^2(P-3Q)}{Q-P}+ \frac{(4\epsilon P^2Q+z^2P+3z^2Q)Py}{z^3(Q-P)}+
\frac{(P+Q)(\epsilon P^2-z^2)^2}{z^6(P-Q)}  \right)dx+\\
\\
  \ \ \ \ \ \ \
  \left( \frac{2QPy^2}{P(Q-P)} + \frac{(3\epsilon P^3+\epsilon P^2Q+z^2P+3z^2Q)Py}{z^3P(P-Q)}+
         \frac{(P+Q)(\epsilon P^2-z^2)^2}{z^6P(Q-P)} \right)dy.\\

\end{array}
\end{equation}
Technical details of compatibility analysis are presented in Appendix.

Since $R=0$, the leaves of the corresponding foliation are special solutions of (\ref{eqgeo3}) and the  dual focal curve  is an arc of the conic $c_0$.
Consider the foliation by integral curves of the vector field $\partial_z+P(z,y)\partial_y$. Its leaves are geodesics (\ref{incidence}) with $D\ne 0$ since these geodesics are transverse to the special solutions $y=const$.  Thus we can normalize $D=1$. As we move in the $z$-direction, the corresponding point $[A:B:C:D]$ in the dual space moves along some curve $\gamma_P$ on the quadric (\ref{quadric}). In the chosen normalization, one easily finds the following parametrization of $\gamma_P$:
\begin{equation}\label{Pgeodesics}
\begin{array}{l}
A=-\frac{1}{P^2}+\frac{\epsilon}{z^2}\\
\\
B= - \frac{z}{P}+\frac{y}{P^2}-\frac{\epsilon y}{z^2}\\
\\
C= -z^2+\frac{2zy}{P}-\frac{y^2}{P^2}+\frac{\epsilon y^2}{z^2}.\\
\end{array}
\end{equation}
Substituting $Q$ for  $P$ in the above expressions, one obtains the family of geodesics $\gamma_Q$, corresponding to the foliations by integral curves of the vector field $\partial_z+Q(z,y)\partial_y$. Easy  exercise in analytic geometry shows that $P,Q,P_y$ satisfies (\ref{3webeq3}) if and only if the arcs $\gamma_P$ and $\gamma_Q$ are planar and belong to the same plane.
\hfill $\Box$\\

The group of projective transformations, generated by vector fields (\ref{pg3}),  naturally acts on the set of all conics of the form (\ref{incidence}) and preserves the quadric (\ref{quadric}). This action   also moves around the planes in $\mathbb R\mathbb P^3$. Namely, the action on the planes
$$
a A+b B + c C+\delta D=0
$$
is generated by the following one-parameter transformation groups $G_i$:
\begin{equation}\label{Planeaction3}
\begin{array}{lllll}
G_1: & a\to a+t_1b+t_1^2c,\ \ \ & b\to b+2t_1^2c,\ \ \ & c\to c, & \delta \to \delta, \\
\\
G_2: & a\to e^{t_2}a, & b\to b, & c\to e^{-t_2}c, & \delta \to \delta, \\
\\
G_3: & a\to a, & b\to b+2t_2a, & c\to c+t_3b+t_3^2a,\ \ \ & \delta \to \delta, \\

\end{array}
\end{equation}
where $t_1,t_2,t_3$ are group parameters.
The action (\ref{Planeaction3}) preservers any quadric $Q_{\mu}$
\begin{equation}
\label{qmu} \mu_1(4ac-b^2)=\mu_2 \delta ^2,
\end{equation}
where $\mu:=[\mu_1:\mu_2]\in \mathbb{RP}^1$. In particular, the (double) plane $\pi_0$, defined by $\delta=0$, is invariant.

\begin{definition}
Two webs are projectively equivalent if there is a projective transformation mapping one to the other. \end{definition}
The hexagonal geodesic 3-webs are parameterized by planes in ${\mathbb R\mathbb P^3}$ different from the plane 
$
\Pi_0=\{[A:B:C:D]\in {\mathbb R\mathbb P^3}: D=0\}. 
$
Therefore the projective moduli space of hexagonal geodesic 3-webs is the space of orbits of points in the complement ${\mathbb R\mathbb P^3}^*\setminus \{\Pi_0\}$, where now $\Pi_0=[0:0:0:1]\in {\mathbb R\mathbb P^3}^*$.
\begin{theorem}
There is a one-parametric family of projectively non-equivalent hexagonal geodesic 3-webs on $(M^2,g)$ with $\dim \mathfrak{p}(g)=3$. Any such hexagonal geodesic 3-web is symmetric with respect to some at least 1-dimensional projective subgroup.
\end{theorem}
{\it Proof:} The orbits of points $[a:b:c:\delta]$ in the compliment ${\mathbb R\mathbb P^3}^* \setminus \{\Pi_0\}$ under the action (\ref{Planeaction3}) are subsets of the quadrics $Q_{\mu}$. On each of these quadrics, the rank of the infinitesimal algebra $\{b\partial _a+2c\partial_b,\ 2a\partial _b+b\partial_c,\ a\partial _a-c\partial_c \}$, generating the action (\ref{Planeaction3}), is equal to 1, if $4ac-b^2=0$, and is equal to 2 at the other points in ${\mathbb R\mathbb P^3}^* \setminus \{\Pi_0\} $. Therefore any connected component of the set $Q_{\mu} \setminus \{[a:b:c:\delta]\in {\mathbb R\mathbb P^3}^*:4ac-b^2=0 \}$ belongs to a single orbit (this is a particular case of Mather's Lemma \cite{Ma-70}). Moreover, any point on that set is stable for some 1-parameter subgroup. Thus the 2-dimensional orbits are parametrized by one parameter $\mu$. 

The 1-dimensional orbits are intersections $Q_{\mu} \cap \{[a:b:c:\delta]\in {\mathbb R\mathbb P^3}^*:4ac-b^2=0 \}$. They also form 1-parameter family, now each point being stable under 2-dimensional subalgebra.   
\hfill $\Box$

\medskip
\noindent {\bf Remark 1}. Like in the case of constant Gaussian curvature, the projective moduli space of hexagonal geodesic 3-webs on $(M^2,g)$ with $\dim \mathfrak{p}(g)=3$ is 1-dimensional, but, in contrast, a generic linear hexagonal 3-web in the Euclidean plane (or generic hexagonal geodesic 3-web on the  surface of constant curvature) does not admit any projective symmetry.

\medskip
\noindent {\bf Remark 2}. To study the geodesic webs, we do not need the metric itself, it is enough to know the form of the equation for geodesics. The following explicit pairwise non-isometric normal forms for pseudo-Riemannian metrics with $\dim \mathfrak{p}(g)=3$ were found in \cite{BMM-08}:
\begin{equation}
\begin{array}{l}
g=\epsilon_1e^{3\bar{x}}d\bar{x}^2+\epsilon_2e^{\bar{x}}d\bar{y}^2,\\
\\
g=\alpha\left(\frac{e^{3\bar{x}}}{(e^{\bar{x}}+\epsilon_2)^2}d\bar{x}^2+\epsilon_1\frac{e^{\bar{x}}}{e^{\bar{x}}+\epsilon_2} d\bar{y}^2\right), \\
\\
g=\alpha\left(\frac{1}{\bar{x}(\gamma \bar{x} +2\bar{x}^2+\epsilon_2)^2}d\bar{x}^2+\frac{\epsilon_1 \bar{x}}{\gamma \bar{x} +2\bar{x}^2+\epsilon_2} d\bar{y}^2\right),
\end{array}
\end{equation}
where $\epsilon_i=\pm 1$,  $\alpha\in \mathbb{R} \setminus \{0\}$, $\gamma \in \mathbb{R}$.

The three operators (\ref{pg3}) belong to three different orbits with respect to the automorphisms of the projective algebra  (there are exactly 3 orbits). 
Any metric with  $\dim \mathfrak{p}(g)=3$ always admits a Killing vector field. The orbit of the Killing vector field  determines which of the three above metrics to take.   The authors of \cite{BMM-08} explain how to recover a metric from the equation for geodesics. This recovery is an integration of some  over-determined system of PDEs, therefore this process brings new constants.

\subsection{Hexagonal geodesic 3-webs on the surfaces with $\dim \mathfrak{p}(g)=2$}

Here  one can choose local coordinates $x,y$ on the surface so that the equation of geodesics is
$$
\frac{d^2y}{dx^2}=H\frac{dy}{dx}+Je^{-2x}\left(\frac{dy}{dx}\right)^3,\ \ \ H,J={\rm const},\ J \ne 0,\ \ H\ne 2,\ H\ne \frac{1}{2}
$$
with the following basis of the projective symmetry algebra:
$$
\{\partial_y,\ \partial_x+y\partial_y \}
$$ (see again \cite{BMM-08}).
Changing the coordinate $x$ for a new coordinate $z$, defined by $\frac{dz}{z}=dx$, and adjusting the integration constant,
we get  the equation of geodesics
\begin{equation}\label{eqgeo2}
z\frac{d^2y}{dz^2}=\rho \frac{dy}{dz}+\epsilon\left(\frac{dy}{dz}\right)^3,\ \  \epsilon=\pm 1, \ \ \rho=const, \ \ \rho\ne 1,\ \ \rho \ne -\frac{1}{2}
\end{equation}
and the symmetry algebra generators
\begin{equation}\label{pg2}
\{\partial_y,\ z\partial_z+y\partial_y\}
\end{equation}
in the coordinates $z,y$. (The excluded values of $\rho$ distinguish the case of  the constant Gaussian curvature.)
Admitting a 2-dimensional symmetry subalgebra, equation (\ref{eqgeo2}) is again integrable in quadratures, but, for most values of $\rho$, not in elementary functions.
Vector fields of the form (\ref{webVfields}), tangent to the web directions, now satisfy the following PDEs:
\begin{equation}\label{2dimPQR}
z(P_z+PP_y)=\rho P+\epsilon P^3,\ \ \ \ \ \ z(Q_z+QQ_y)=\rho Q+\epsilon Q^3,\ \ \ \ \ \ z(R_z+RR_y)=\rho R+\epsilon R^3.
\end{equation}
The geometry of hexagonal geodesic 3-webs for  $\dim \mathfrak{p}(g)=2$ is not as rich as that one for  $\dim \mathfrak{p}(g)=3$.
\begin{theorem}\label{webdim2}
Suppose that $\dim \mathfrak{p}(g)=2$ and the coordinates $(z,y)$ on the surface $(M^2,g)$ are chosen as above. Then a geodesic 3-web ${\cal{G}}_3$ on $M^2$ is hexagonal if and only if it admits an infinitesimal symmetry $\partial_y,$ and the vector fields (\ref{2dimPQR}) (renamed, if necessary) satisfy $R=P+Q=0$.
\end{theorem}
{\it Proof:} The Blaschke curvature
 of our web vanishes if
\begin{equation}\label{cur02}
Pyy+Qyy+Ryy-\frac{3\epsilon}{z}(PPy+QQy+RRy)+\frac{1}{z^2} [P^3+Q^3+R^3+\epsilon (\rho-1)(P+R+Q)]=0.
\end{equation}
Analysis of compatibility of the constraint (\ref{cur02}) with system (\ref{2dimPQR}) (see Appendix for more detailed discussion) implies that:\\
 1) one of the slopes $P,Q,R$ vanishes identically,\\
 2) the other two sum up to zero,\\
 3) the slopes do not depend on $y$.\\
  We may suppose that $R=0$, $P+Q=0$. Then our hexagonal geodesic 3-web is completely determined by $P$, satisfying
 $P_y=0$, $zP_z=\rho P+\epsilon P^3$. (Note that this system is easily integrated in elementary functions, but we do not need explicit formulas.)
\hfill $\Box$\\

\noindent Now let us describe symmetry properties of the webs under consideration.
\begin{theorem}
The projective moduli space of hexagonal geodesic 3-webs on $(M^2,g)$ with $\dim \mathfrak{p}(g)=2$ is discrete. 
For $\epsilon \rho < 0$, there is a unique hexagonal geodesic 3-web stable under the action of the whole 2-dimensional projective algebra, the slopes of this web being constant.
\end{theorem}
{\it Proof:} Observe that the infinitesimal symmetry $z\partial_z+y\partial_y$ preserves the slopes $P,Q,R$ and changes $z$, therefore the projective moduli space of hexagonal geodesic 3-webs is discrete. If such a web is symmetric with respect to this symmetry then the slopes are constant. This web exists only for $\epsilon \rho < 0$ since the non-zero slopes $P,Q=-P$ are the roots of the quadratic equation $\rho+\epsilon P^2=0.$
\hfill $\Box$

\medskip
\noindent {\bf Remark}. The following explicit pairwise non-isometric normal forms for pseudo-Riemannian metrics with $\dim \mathfrak{p}(g)=2$ were found in \cite{BMM-08}:
\begin{equation}
\begin{array}{l}
g=\epsilon_1e^{(\beta+2)\bar{x}}d\bar{x}^2+\epsilon_2e^{\beta \bar{x}}d\bar{y}^2,\\
\\
g=\alpha\left(\frac{e^{(\beta+2)\bar{x}}}{(e^{\beta\bar{x}}+\epsilon_2)^2}d\bar{x}^2+\epsilon_1\frac{e^{\beta\bar{x}}}{e^{\beta\bar{x}}+\epsilon_2} d\bar{y}^2\right), \\
\\
g=\alpha\left(\frac{e^{2\bar{x}}}{\bar{x}^2}d\bar{x}^2+\frac{\epsilon}{\bar{x}} d\bar{y}^2\right),
\end{array}
\end{equation}
where $\epsilon,\epsilon_i=\pm 1$, $\beta\in \mathbb{R} \setminus \{-2,0,1\}$, $\alpha\in \mathbb{R} \setminus \{0\}$. The choice of the form depends, in particular, on the parameter $\rho $ in equation (\ref{eqgeo2}).
The authors of \cite{BMM-08} explain how to recover a metric from the equation for geodesics (\ref{eqgeo2}). This recovery is an integration of some  over-determined system of PDEs therefore this process brings a new constant $\alpha$. Essential is that $\bar{x}$  depends linearly only on $x$ and $\bar{y}$ depends linearly only on $y$. Thus each hexagonal geodesic 3-web is symmetric, in fact, with respect to the Killing vector field $\partial_y$.

\section{Old and new examples} 

For a given metric, the space $\mathcal{I}_3$ of cubic integrals is finite-dimensional (see \cite{Kr-08}), the maximal possible dimension being 10. This bound is sharp: the family of hexagonal geodesic 3-webs is 9-parametric on the surfaces of constant Gaussian curvature, as was known already a century ago (see \cite{Ma-31}). The case of metrics of non-constant curvature is less understood: the bound of 7  for the dimension of $\mathcal{I}_3$ was proved and the sharp bound of 4 was conjectured in \cite{Kr-08}. In fact, all known  examples of surfaces with ${\rm dim}( \mathcal{I}_3)>1$ or carrying parametric families of hexagonal geodesic 3-webs are 
the surfaces of constant Gaussian curvature (9-parameter families)
and the surfaces admitting Killing vector field (1-parameter family \cite{Fi-99},  ${\rm dim}( \mathcal{I}_3)=4$ \cite{BMM-08,MS-11,VDS-15}), in particular, the surfaces with 3-dimensional projective group, also known as superintegrable (${\rm dim}( \mathcal{I}_3)=4$).  

In this section we review the known examples and present some new ones, a special attention will be paid  to the surfaces carrying many hexagonal geodesic 3-webs.    

As it was already mentioned, the general solution to (\ref{hydro}) can be found by generalized hodograph method, discovered by Tsarev \cite{Ts-85}. This method needs an integration of  overdetermined (but consistent!) linear system of PDEs for commuting flows, the flows being determined by functions, depending only on Riemann invariants. We will not pursue this line here. Instead we discuss how to obtain some classes of solutions  by imposing additional constraints, compatible with (\ref{hydro}). 

\subsection{Simple waves}
Existence of Riemann invariants implies compatibility of very simple finite (not differential) constraints. Namely, it is immediate from the form (\ref{n}) that setting one or two Riemann invariants constant, one gets a consistent system. Note that if all 3 Riemann invariants are constant, then the metric is flat. Setting exactly one of the Riemann invariants constant reduces the number of field variables of  system (\ref{hydro}) to two. Such systems are linearized via classical hodograph transform (see \cite{RJ-83}).   

With $R^1=c_1$, $R^2=c_2$ one reduces (\ref{n}) to a single Hopf equation for $R=R^3$: 
$$
R_u=\lambda (R)R_v,
$$
where $\lambda(R)=\lambda^3(c_1,c_2,R)$. The well known general solution depends on one arbitrary function $f(R)$ and reads as 
\begin{equation}\label{solHopf}
\lambda(R)u+v=f(R). 
\end{equation}
Thus we can define a local solution by the following procedure:\\
1) solve equation (\ref{lambdas}) for $\lambda^1$,\\
2) express Riemann invariants via $F,\lambda^2,\lambda^3$, (see (\ref{R3})),\\
3) find $F,\lambda^2$ as functions of $\lambda^3$ from the constraints $R_1=c_1$, $R_2=c_2$,\\
4)  get $R=R^3$ as a function of $\lambda^3$,\\
5) substitute the found expressions to (\ref{solHopf}) and find $\lambda^3$ as a function of $u,v$,\\
6) use Vieta's formulas for the characteristic equation (\ref{char}) to get $E$ and $G$. $F$ is already found on the step 3).\\
Choosing $f$, one have to take care that $E>0,\  EG-F^2>0$. Some steps in the described procedure cannot be made explicitly, therefore solutions  are effectively defined implicitly by two equations: one is (\ref{solHopf}) and the other relates $\lambda^2$ and $\lambda^3$ after excluding $F$ between the equations $R_1=c_1$, $R_2=c_2$.
\smallskip 

\noindent {\bf Remark.} Any solution to (\ref{hydro}) parametrizes some (possibly singular) manifold in the hodograph space: 
$$
(u,v)\mapsto (E(u,v),F(u,v),G(u,v)).
$$ 
Generically, the manifold is a  surface, but can degenerate to a curve. Such solutions with degenerate hodograph are called simple waves (see \cite{RJ-83}).  The solutions described above have 1-dimensional hodograph, which is a straight line if one chooses the Riemann invariants as the coordinates in the hodograph space. Although some specific simple wave solutions give metrics with Killing vector filed (see the next subsection), a general solution of the type does not manifest non-trivial symmetry properties.

\subsection{Group invariant solutions}
Let us describe solutions invariant with respect to one-dimensional symmetry subalgebras. According to 
Proposition \ref{Liesysm}, it is enough to consider generators in one of the forms (\ref{normalLiesysm}). In fact, solutions invariant with respect to any other symmetry subalgebra can be obtained by a convenient transformation of the whole symmetry group (see \cite{Ov-82}). Among all the solutions we should select those that satisfy obvious geometric restrictions. Namely, the quadratic form $E(u,v)du^2+2F(u,v)dudv+G(u,v)dv^2$ must be positively definite and non-degenerate. Moreover, we will be interested mostly in metrics with non-constant Gaussian curvature. \\ 

\noindent$\bullet$ {\it Subalgebra generated by operator $\partial_u+\kappa \partial_v$ }\\

For an invariant solution, the operator $\partial_u+\kappa \partial_v$ is a Killing vector field on the corresponding surface. A hexagonal geodesic web, invariant with respect to this vector field flow, is one from at least one-parametric family of such webs: in fact, it is well known that surfaces admitting a Killing vector field can be locally realized as surfaces of revolution, which always have one-parametric rotationally invariant family of hexagonal geodesic 3-webs (see \cite{Fi-99}).  Note that the corresponding solutions have one-dimensional hodograph and therefore they are particular cases of simple waves from the preceding subsection.

Invariant solutions have the form $E=e(s),\ G=j(s),\ F=f(s),$ where $s=v-\kappa u$ is the invariant of the subalgebra. Substituting this Ansatz into (\ref{hydro}), one obtains a linear homogeneous system  for the derivatives 
of $e,j,f$. For a solution to define a surface with non-flat metric, the discriminant of this system must vanish. This condition is equivalent to  
$$
(2\kappa +1)e+\kappa (\kappa -1)f-\kappa ^2(\kappa +2)j=0.
$$ 
Elementary analysis of compatibility of this equation with geometric restrictions shows that non-constant solutions exist only for 
the following values of the parameter: $\kappa =1$, $\kappa =-\frac{1}{2}$, $\kappa =-2$. 
Observe that these three values belong to the same orbit of the discrete symmetry group (see Section \ref{symsec}). Therefore one can choose $\kappa =1$.   
\begin{proposition}\label{translationpur}
Any (local) solution, invariant with respect to the operator $\partial_u+\partial_v$, has the form 
$$
E=h^2(s),\ G=h^2(s), \ F=f_0 h(s)-h^2(s)
$$
with $s=v- u$,  constant $f_0\ne 0$ and an arbitrary function $h$. 
The Gaussian curvature $K_G$ of the corresponding surface is  
$$
K_G=\frac{h''}{h^2(f_0-2h)}+\frac{(3h-f_0)(h')^2}{h^3(f_0-2h)^2}.
$$

\end{proposition}
{\it Proof:} For $\kappa =1$ one obtains $e=j$ from the condition on the discriminant and the following ODE $2j\frac{df}{ds}=(f-j)\frac{dj}{ds} $ from (\ref{hydro}), which is easy to integrate.  
\hfill $\Box$\\

\noindent {\bf Remark 1.} Metrics, admitting a Killing vector and more then 2 independent cubic integrals of its geodesic flow, were classified in \cite{MS-11}: it turns out that the space $\mathcal{I}_3$ of cubic integrals is of dimension 4. The classification was obtained by reducing the overdetermined system of PDEs for the integrals to the system of ODEs. The reduction is based on the following observation: the Killing vector field gives an integral $I_1$, linear in momentum (Noether's Theorem), this integral defines a linear action on the finite-dimensional space $\mathcal{I}_3$ by $I_3\mapsto \{I_1,I_3\}_H$. There is 2-dimensional eigensubspace of $\mathcal{I}_3$ generated by $I_1H,I_1^3$, the eigenvalue being zero. Thus the existence of additional integrals implies the existence either of eigenvalues different from zero, or the existence of an integral $I_3$ satisfying $\{I_1,I_3\}_H=c_1I_1H+c_2I_1^3$ with some constant $c_1,c_2$. 
This trick generalizes straightforwardly to the case of any infinitesimal symmetry of system (\ref{hydro}). In fact, any symmetry of (\ref{hydro}) is a symmetry of the geodesic flow by Theorem \ref{projsym} and therefore maps integrals to integrals. Moreover, this symmetry acts linearly on $\mathcal{I}_3$. Unfortunately, the complete analysis is relatively simple only for the case of symmetry in form (\ref{normalLiesysm}.1)) and rather involved for the other two forms of symmetry. We do not go through all the details and  present here only some  
partial results of "symbolic experiments" with Maple. There are  doubts that this analysis would brings new interesting forms: we have only one free parameter, namely $\kappa$, and initial values of some ODEs to play with, whereas for the case (\ref{normalLiesysm}.1)), an arbitrary function of one variable was in our disposal.

\medskip

\noindent$\bullet$ {\it Subalgebra generated by operator $\partial_u+\kappa \partial_v+E\partial_E+F\partial_F+G\partial_G$ }\\   

Invariant solutions have the form $E=\exp(u)e(s),\ G=\exp(u)j(s),\ F=\exp(u)f(s),$ where $s=v-\kappa u$. Substituting this Ansatz into (\ref{hydro}), one obtains a linear system  for the derivatives 
of $e,j,f$. If the discriminant of this system does not vanish, which is equivalent to  
$$
\delta :=-j\kappa^3+(f-2j)\kappa^2+(2e-f)\kappa+h\neq 0,
$$ then one obtains 
\begin{equation}\label{sym2}
\begin{array}{l}
e'=\frac{1}{\delta}(\kappa+1)e(f-\kappa j),\\
\\
j'=\frac{1}{\delta}j(-j\kappa^2+(2f-2j)\kappa+e),\\
\\
f'=\frac{1}{2\delta}(-2fj\kappa^2+(2f^2-3fj+je)\kappa+e(f+j)).

\end{array}
\end{equation} 
Elementary analysis shows that the condition $\delta=0$ is not compatible with the geometric restrictions. 
\begin{proposition}\label{translationcentre}
Any (local) solution to (\ref{hydro}), invariant with respect to the operator $\partial_u+\kappa \partial_v+E\partial_E+F\partial_F+G\partial_G$, is given by $E=\exp(u)e(s),\ G=\exp(u)j(s),\ F=\exp(u)f(s),$
where $s=v-\kappa u$ and  $e,j,f$ solve ODEs (\ref{sym2}). The Gaussian curvature of the corresponding surface is constant if and only if $\kappa \in \{0,-1\}$. Moreover, for these values of $\kappa $, the metric is flat.   
\end{proposition}
{\it Proof:} The Gaussian curvature $K_G$ is 
$$
K_G=\frac{\kappa (\kappa+1)[j(2f-j)\kappa^4+2j(2f-j)\kappa^3-2e(2f-e)\kappa-e(2f-e)]}{4\exp(u)\delta ^3}.
$$ 
Therefore the curvature is constant if and only if it vanishes.
Direct calculation shows that the expression in the square brackets can not be zero for "geometric" solutions.
\hfill $\Box$\\

\noindent {\bf Remark 2.} By Lemma \ref{spiralLie},  hexagonal geodesic 3-webs with the symmetry subalgebra,  considered in Propostion \ref{translationcentre}, can be realized on Lie  spiral surfaces. Hexagonal geodesics 3-webs on  Lie  spiral surfaces were studied in \cite{Sa-26}. Unfortunately, the treatment was incomplete: the author claims that this type of symmetry is impossible.   

\begin{proposition}\label{translationcentrek1}
Any (local) solution to (\ref{hydro}), invariant with respect to the operator $X=\partial_u+ \partial_v+E\partial_E+F\partial_F+G\partial_G$, determines a metric  admitting 3-dimensional projective algebra. 
\end{proposition}
{\it Proof:} Due to the classification results of \cite{MS-11}, it is enough to show that the space $\mathcal{I}_3$ is 4-dimensional and any cubic integral is divisible by a linear one. 

Let us look for the solution in the following form: $E=\exp\left(\frac{u+v}{2}\right)e(s),\ G=\exp\left(\frac{u+v}{2}\right)j(s),\ F=\exp\left(\frac{u+v}{2}\right)f(s),$ where $s=v-u$. Then
$$
e'=\frac{(4f-j-3e)e}{6(e-j)},\ \ j'=\frac{(4f-3j-e)j}{6(e-j)}, \ \ 
f'=\frac{(f-e)(f-j)}{3(e-j)}.
$$
If $\xi \partial_u+\eta \partial _v$ is a tangent vector, then the action of $X$ on $\xi, \eta$ is trivial.  
The relations
\begin{equation}\label{velocitymomentum}
E\xi+F\eta=p, \ \ F\xi+G\eta=q,
\end{equation}
give the action of $X$ on the momentum: $X(p)=p$, $X(q)=q$,
and therefore on cubic integrals 
$$
X(I_3)=X(K_3p^3+K_2p^2q+K_1pq^2+K_0q^3)=X(K_3)p^3+X(K_2)p^2q+X(K_1)pq^2+X(K_0)q^3)+3I_3.
$$
In particular, Lemma (\ref{LMcub}) allows to compute the action $X(I)=2I $ on the integral
(\ref{webcubicint}). Substituting the Ansatz $K_i=\exp\left(-\frac{u+v}{2}\right)k_i(s)$, $i=0,1,2,3,$ into PDEs for cubic integrals, one obtains 5 ODEs for four derivatives $k'_i$. Compatibility analysis shows that the space of solutions is 2-dimensional. Therefore the eigensubspace of $\mathcal{I}_3$ corresponding to the eigenvalue $\lambda_1=2$ is 2-dimensional.   

In a similar way, one checks that there is one-dimensional space of linear integrals $I_1=
\exp\left(-\frac{u+v}{6}\right)(\alpha (s)p+\beta (s)q)$, one-dimensional eigensubspace of cubic integrals corresponding to the eigenvalue $\lambda_2 =\frac{7}{3}$, and one-dimensional eigensubspace of cubic integrals corresponding to the eigenvalue $\lambda_3 =\frac{5}{3}$. Finally, one sees from the obtained expressions for cubic integrals, that all of them are divisible by $I_1$. 
\hfill $\Box$\\

\noindent$\bullet$ {\it Subalgebra generated by operator $u\partial_u+v\partial_v+\kappa(E\partial_E+F\partial_F+G \partial_G)$ }\\  

 For $\kappa \ne -2$ the surfaces, corresponding to the solutions invariant with respect to this operator, can be immersed into space as Lie spiral surfaces  (Lemma \ref{spiralLie}). Such surfaces, carrying hexagonal geodesic 3-webs, were also studied by Sauer \cite{Sa-26}.

Invariant solutions have the form $E=u^{\kappa}e(s),\ G=u^{\kappa}j(s),\ F=u^{\kappa}f(s),$ where $s=\frac{v}{u}$. Substituting this Ansatz into (\ref{hydro}), one obtains a linear system  for the derivatives 
of $e,j,f$. If the discriminant of this system do not vanish, which is equivalent to   
$\Delta: =-js^3+(f-2j)s^2+(2e-f)s+e \neq 0$, then one obtains 
\begin{equation}\label{sym3}
\begin{array}{l}
e'=\frac{1}{\Delta}\kappa(s+1)e(f-js),\\
\\
j'=\frac{1}{\Delta} \kappa j (-js^2+(2f-2j)s+e),\\
\\
f'=\frac{1}{2\Delta}\kappa(-2fjs^2+[2f^2-3fj+je]s+e[f+j]).

\end{array}
\end{equation} 
\begin{proposition}\label{dilataion}
Any (local) solution to (\ref{hydro}), invariant with respect to the operator $u\partial_u+v\partial_v+\kappa(E\partial_E+F\partial_F+G\partial_G)$, is given by $E=u^{\kappa}e(s),\ G=u^{\kappa}j(s),\ F=u^{\kappa}f(s),$ where $s=\frac{v}{u}$ and the functions $e,j,f$ are defined as follows.

If $\Delta \equiv 0$ then $\kappa =0$, the functions $e,j,f$ satisfy 
\begin{equation}\label{sym3discr}
\begin{array}{l}
e=\frac{s(js^2+(2j-f)s+f)}{2s+1},\\
\\
j'=\frac{2j(js+2j-3f)(4js^3+(7j-2f)s^2+(4j-2f)s+f)}{(s-1)(2s+1)(s+2)(js+f)^2},\\
\\
f'=\frac{(f-js)(js+2j-3f)(4js^3+(7j-2f)s^2+(4j-2f)s+f}{(s-1)(2s+1)(s+2)(js+f)^2},

\end{array} 
\end{equation} 
and the Gaussian curvature of the corresponding surface is not constant. 

If $\Delta \ne 0$ and $\kappa \ne -2$ then the functions $e,j,f$ solve ODE (\ref{sym3}) and the Gaussian curvature of the corresponding surface is constant, namely zero, only for $\kappa =0$.   

If $\Delta \ne 0$ and $\kappa =-2$, the Gaussian curvature is constant if and only if the solution of (\ref{sym3}) is subject to one of the following constraints, which are compatible with (\ref{sym3}): 
\begin{equation}
\begin{array}{llcl}
a) & (2f-j)s^2+(2s+1)e=0, & with & K_G=-\frac{j}{(js+f)^2},\\
\\
b) & (s^2+2s)j+2f-e=0, & with & K_G=-\frac{e}{(s+2)^2(js+f)^2},\\
\\
c) & js^2-e=0, & with & K_G=\frac{2f-e-j}{(s-1)^2(js+f)^2}.
\end{array}
\end{equation} 
\end{proposition}
{\it Proof:} The case $\Delta = 0$ is straightforward: one resolves this equation for $e$ and obtains $\kappa =0$ as the compatibility condition.  

The analysis of the Gaussian curvature for general case becomes rather involved and we use computer algebra system Maple to carry out the necessary calculations. The Gaussian curvature $K_G$ has the form $K_g=\frac{1}{\Delta ^3}\kappa u^{\kappa+2}N(e,f,j,s,\kappa)$, where $N$ is polynomial. Thus if the curvature is constant for $\kappa \notin \{-2,0\}$ than it vanishes. 
Differentiating the equation $N=0$ two times and excluding $e,f$ between thus obtained three equations, one 
obtains that necessarily $\kappa \in \{-6,-5,-3,-2,0,1,3\}$. The case study shows that  the values different from $-2$ are not compatible with geometric restrictions. 

For $\kappa =-2$ the numerator of equation $\frac{dK_G}{ds}=0$ factors into four polynomials, only three 
of them being compatible with system (\ref{sym3}) and geometric restrictions.
\hfill $\Box$\\

 For $\kappa =-2$ the operator gives a Killing vector field on the surface that can be realized locally as a surface of revolution.  The hexagonal geodesic 3-webs on such surfaces, invariant with respect to this operator, were described first in \cite{Sa-26}. The invariant 3-web is not symmetric with respect to meridians and therefore not of the type described by Finsterwalder \cite{Fi-99}. We complete the results of \cite{Sa-26} and show that the family of webs on any of the described surfaces is, in fact, 3-parametric.      
 
\begin{proposition}
Any (local) solution to (\ref{hydro}), invariant with respect to the operator $X=u\partial_u+v\partial_v-2(E\partial_E+F\partial_F+G \partial_G)$ determines a metric  admitting a Killing vector field. The geodesic flow of the metric has a 4-dimensional space of cubic integrals. 
\end{proposition}
{\it Proof:} We apply again the trick indicated in Remark 1. 
Consider the solution as in Proposition \ref{dilataion}.  
The action $X(\xi)=\xi$, $X(\eta)=\eta$ on the coordinates $\xi,\eta$ of the tangent vector $\xi \partial_u+\eta \partial _v$ determines again the action on the momentum: $X(p)=-p$, $X(q)=-q$ by formulas (\ref{velocitymomentum}).
For cubic integrals we have 
$$
X(I_3)=X(K_3p^3+K_2p^2q+K_1pq^2+K_0q^3)=X(K_3)p^3+X(K_2)p^2q+X(K_1)pq^2+X(K_0)q^3)-3I_3.
$$
 Substituting the Ansatz $K_i=\exp(\nu u)k_i(s)$, $i=0,1,2,3,$ where $s=v/u$, into PDEs for cubic integrals and studying the compatibility conditions,  one checks that there is a two-dimensional eigensubspace of $\mathcal{I}_3$ corresponding to the eigenvalue $\lambda_1=0$, a one-dimensional eigensubspace for $\lambda_2=1$,   and a one-dimensional eigensubspace for $\lambda_3=-1$. 
Thus a generic surface is of the type described in \cite{MS-11} and its geodesic flow admits a 4-dimensional space of cubic integrals.
\hfill $\Box$\\

Finally, let us show that invariance with respect to 2-dimensional subalgebras does not bring new interesting solutions. 
\begin{proposition} 
Any (local) solution to (\ref{hydro}), invariant with respect to some two-dimensional symmetry subalgebra, defines a surface with constant Gaussian curvature.
\end{proposition}  
{\it Proof:} Any 2-dimensional Lie algebra is either commutative or possesses  a basis $X,Y$ such that holds true $[X,Y]=X$. Thus, for the operator $X$, we can choose one considered in Propositions \ref{translationpur}, \ref{translationcentre}, \ref{dilataion}.\\
\noindent$\bullet$ {\it Case $X=\partial_u+\partial_v$ }\\ 
If $[X,Y]=0$ then the operator $Y$ has the form $Y=c_1T_1+c_2T_2+c_4D_2$ (see formula (\ref{sym}) for the notation), where at least one of the coefficients $c_1,c_2$ does not vanish. Subtracting from $Y$ a suitable multiple of $X$ and transposing $u$ and $v$, if necessary, we can bring $Y$ to the form  $Y=c_2T_2+c_4D_2$, where still $c_2\ne 0$. The coefficient $c_4$ is also non-vanishing: otherwise the surface admits two Killing vector fields and therefore has constant Gaussian curvature. Now applying $\exp_{tD_1}$ (and changing the signs of $u$ and $v$, if necessary) we get $Y$ in the form  $ \partial_v+E\partial_E+F\partial_F+G\partial_G.$   
Solutions, invariant with respect to this 2-dimensional algebra, can be written as
$$
E=e_0\exp(v-u),\ G=j_0\exp(v-u), \ F=f_0\exp(v-u)
$$ 
with some constants $e_0>0,j_0>0,f_0$, satisfying $e_0j_0-f_0^2>0$. But equations (\ref{hydro}) give $e_0=j_0=-f_0$ violating the non-degeneracy condition.\\
If $[X,Y]=X$ then the operator $Y$ has the form $Y=c_1T_1+c_2T_2+D_1+c_4D_2$. As in the commutative case, applying linear substitutions for $Y$ and action of the whole symmetry group, one brings $Y$ to the form 
$u\partial_u+ v\partial_v+\kappa (E\partial_E+F\partial_F+G\partial_G).$ Substituting the following Ansatz for the invariant solutions
$$
E=e_0(v-u)^{\kappa},\ G=j_0(v-u)^{\kappa}, \ F=f_0(v-u)^{\kappa}
$$    
into (\ref{hydro}), one argues that either $\kappa =0$ (and the metric is flat) or $e_0=j_0=-f_0$ (and the searched quadratic form is degenerate).

\noindent$\bullet$ {\it Case $X=\partial_u+\kappa \partial_v+E\partial_E+F\partial_F+G\partial_G$}\\
If $[X,Y]=0$ then the operator $Y$ has again the form $Y=c_1T_1+c_2T_2+c_4D_2$, where at least one of the coefficients $c_1,c_2$ do not vanish. Substituting $Y$ by $Y-c_4X$ we arrive at the case already considered above. 
The commuting relation $[X,Y]=X$ is not possible since the derivative of the whole symmetry algebra is spanned by $T_1,T_2$.\\
\noindent$\bullet$ {\it Case $X=u\partial_u+v\partial_v+\kappa(E\partial_E+F\partial_F+G\partial_G)$}\\
If $[X,Y]=0$ then the operator $Y$ has the form $Y=c_3D_1+c_4D_2$ with  $c_3\ne 0$. Substituting $Y$ by $Y-c_3X$ we get a symmetry operator proportional to $D_2$, which can not be tangent to graphs of solutions. 
The commuting relation $[X,Y]=X$ is again not possible for this $X$. 
\hfill $\Box$

\section{Concluding remarks}

\vspace{1pt} \noindent$\bullet$ {\it Polynomial integrals of higher degrees and in higher dimensions}

\smallskip
It seems that Theorem \ref{cubandflat} can not be straightforwardly generalized for polynomial integrals of higher degree. Generically, the webs defined by these integrals do not have non-trivial Abelian relations and, consequently, are not parallelizable. 

It would be also interesting to understand which of the results of this paper survive in higher dimensions.  

\smallskip

\vspace{1pt} \noindent$\bullet$ {\it Distinguished coordinates} 

\smallskip
 
If one chooses special coordinates on the surface, such as semi-geodesic or Chebyshev (see \cite{BM-11,MP-17}), or isothermic (see \cite{Ko-82}), then the quasilinear system, describing a metric and a cubic first integral of the corresponding geodesic flow, involve 3 field variables and is diagonalizable,  semi-Hamiltonian and therefore integrable. For cubic integrals, Theorem \ref{Tsarevintegrability} extends this list to the "flat web" coordinates. (In these coordinates  $u,v$, the 3-web foliations are $u=const$, $v=const$, $u+v=const$.)      

\smallskip

\noindent$\bullet$ {\it Pseudo-Riemannian metrics}
\smallskip

The correspondence between cubic integrals of geodesic flow and hexagonal geodesic 3-webs on the surface, as well as the established properties of the PDEs for the metric,  hold true for pseudo-Riemannian metrics.

\smallskip

\section*{Acknowledgements}
The author thanks Ferapontov E.V. and Matveev V.S. for fruitful discussions.
This research was supported by grant \#2017/02954-2 of S\~ao Paulo Research Foundation (FAPESP).
The author thanks also the hospitality of the Laboratoire de Mathématiques de Reims, where this study was accomplished, and V. Ovsienko in particular.

\section{Appendix}
Here we present the details of computations, mentioned in the proofs.

\subsection{Formulas for the proof of Theorem \ref{cubandflat}}
The 5 equations, obtained from the commuting condition $\{I,H\}=0$, can be resolved for the following derivatives:
\begin{equation}
\begin{array}{l}
E_u=\frac{E}{5K}   \left(\frac{2(KF-LG)}{EG-F^2}E_v+\frac{4(LF-KE)}{EG-F^2}F_v +\frac{2KEF+3LF^2-5LEG}{G(EG-F^2)}G_v-2L_v-2K_u\right),\\
\\
F_u= \frac{F}{5K}\left(\frac{5KEG-2LFG-3KF^2}{2F(EG-F^2)}E_v +\frac{2(LF-KE)}{EG-F^2} F_v+\frac{2KEF-5LEG+3LF^2}{2G(EG-F^2)}G_v-L_v-K_u\right),\\
\\
G_u= 2F_v-\frac{F}{G}G_v,\\
\\
L_u=\frac{2L}{5K}\left( \frac{2LG+3KF}{EG-F^2}E_v -\frac{2(2LF+3KE)}{EG-F^2}F_v +\frac{3KEF+5LEG-3LF^2}{G(EG-F^2)}G_v+2L_v+2K_u\right),\\
\\
K_v=\frac{2K}{EG-F^2}\left(2FF_v -GE_v-EG_v\right).
\end{array}
\end{equation}
Differentiating the above formulas, one computes all the mixed derivatives of $E,F,G,L,K$, the second order derivatives of $E,F,G,L$ with respect to $u$ and the derivative $K_{vv}$ in terms of $E,F,G,L,K,E_v,F_v,G_v,L_v,K_u,E_{vv},F_{vv},G_{vv},L_{vv},K_{uu}$. Substituting the found expressions into 
$\frac{\partial}{\partial v}\left(\frac{K_u}{K}\right)=\frac{\partial}{\partial u}\left(\frac{L_v}{L}\right)$ we obtain an identity. 

\subsection{Compatibility conditions in the proof of Theorem \ref{analogGS}}

Equations (\ref{3dimPQR}) allows one to find all mixed derivatives and derivatives with respect to $z$ via derivatives with respect to $y$. Equation (\ref{cur03}) gives $P_{yy}$. One  has $d(dP_y)=0$, which gives $Q_{yyy}$. Similarly, $d(dQ_{yy})=0$ gives $R_{yyyy}$. Finally,  $d(dR_{yyy})=0$ gives a constraint of unexpectedly low, namely, second order. We resolve this constraint for $R_{yy}$ and compare the $z$-derivative of thus obtained expression for $R_{yy}$ with already known expression for $R_{zyy}$, which gives $PR=0$. Thus, we can assume $R=R_y=0$. Now the expression for $R_{yy}$ depends on derivatives of first order. Resolving the equation $R_{yy}=0$ for $Q_y$, we arrive at system (\ref{3webeq3}).

\subsection{Compatibility conditions in the proof of Theorem \ref{webdim2}}

Equations (\ref{2dimPQR}) allows one to find all mixed derivatives and derivatives with respect to $z$ via derivatives with respect to $y$. Equation (\ref{cur02}) gives $P_{yy}$. One  has $d(dP_y)=0$, which gives $Q_{yyy}$. Similarly, $d(dQ_{yy})=0$ gives $R_{yyyy}$. Finally,  $d(dR_{yyy})=0$ gives a constraint of unexpectedly low, namely, second order. Resolving this constraint for $R_{yy}$ and comparing the $z$-derivative of thus obtained expression for $R_{yy}$ with already known expression for $R_{zyy}$, we find $Q_y$. Comparing the $z$-derivative of thus obtained expression for $Q_y$ with already known expression for $Q_{zy}$, we get $R$. Comparing the $z$-derivative of thus obtained expression for $R$ with already known expression for $R_z$, we finally obtain $$
(\rho-3)(\rho-1)(3\rho-1)(2\rho+1)(5\rho+3)(\rho+2)(\rho+3)QP(Q+P)=0.
$$
The values $1,-\frac{1}{2}$ of $\rho$ correspond to the case of constant Gaussian curvature. The above described scheme of analysis, applied a pair of steps further for the other values, leads to the same conclusion: one of the slope is zero, the other two sum up to zero. Finally, one checks that one of these conditions implies the other. Moreover, none of the slopes depends on $y$.

In the presented compatibility analysis, one also has to check that the vanishing of denominators of the expressions for derivatives does not give any new case. Checking this, it is useful to keep in mind the following simple observation: if two of the slopes are constant, then Theorem \ref{webdim2} holds true.

\end{document}